\documentclass{amsart}
\newcommand{\contains}{\supset}

\newcommand{\bbP}{{\mathbb P}}
\newcommand{\bbQ}{{\mathbb Q}}
\newcommand{\bbC}{{\mathbb C}}
\newcommand{\Qbar}{\overline{{\mathbb Q}}}
\newcommand{\calC}{{\mathcal C}}
\newcommand{\calD}{{\mathcal D}}
\newcommand{\calE}{{\mathcal E}}
\newcommand{\calS}{{\mathcal S}}
\newcommand{\calL}{{\mathcal L}}
\newcommand{\Oh}{{\mathcal O}}
\newcommand{\cohom}{{\mathrm H}}
\newcommand{\abs}[1]{{\mid\,#1\,\mid}}
\theoremstyle{plain}
\newtheorem*{extthm*}{Theorem}
\newtheorem{thm}{Theorem}
\newtheorem{lemma}[thm]{Lemma}
\newtheorem{propn}[thm]{Proposition}
\newtheorem{cor}[thm]{Corollary}

\theoremstyle{definition}

\newtheorem{rmk}{Remark}
\newtheorem{prob}{Problem}

\newtheorem*{ack*}{Acknowledgments}
\newcommand{\ie}{i.~e.\ }
\begin{document}
\title{A Geometric characterization of Arithmetic Varieties}
\author{Kapil Hari Paranjape}
\begin{abstract}
  A result of Belyi can be stated as follows. Every curve defined over
  a number field can be expressed as a cover of the projective line
  with branch locus contained in a rigid divisor. We define the notion
  of geometrically rigid divisors in surfaces and then show that every
  surface defined over a number field can be expressed as a cover of
  the projective plane with branch locus contained in a geometrically
  rigid divisor in the plane. The main result is the characterisation
  of arithmetically defined divisors in the plane as geometrically
  rigid divisors in the plane.
\end{abstract}
\maketitle
\section{Introduction}
This paper is an attempt to generalise a result of
Belyi (see \cite{Belyi}).
\begin{extthm*}[Belyi]
  Let $C$ be a smooth projective curve over an algebraic number field
  and $T$ a finite set of closed points in $C$. There is a finite
  morphism $f:C\to\bbP^1$ so that the image $f(T)$ and the branch
  locus of $f$ are contained in the set of three points
  $\{0,1,\infty\}$.
\end{extthm*}
We note that this gives a completely geometric characterisation of
algebraic curves over number fields, since any deformation of a triple
of points in $\bbP^1$ is in fact trivialised by an automorphism of
$\bbP^1$.

A naive generalisation of this could require a surface over a number
field to be expressible as a cover of $\bbP^2$ that is \'etale outside
four general lines; however, as Koll\'ar pointed out, this fails since
the fundamental group of the complement of four general lines in
$\bbP^2$ is abelian, whereas many surfaces have non-abelian
fundamental groups.  Thus one needs to look at more general divisors
in $\bbP^2$. The problem is that these divisors have non-trivial flat
deformations. We need to find an algebraic notion that restricts the
possible deformations. Thus, in Section 1 we define the notion of {\em
  geometrically rigid} divisors on a surface.

Let $C$ be any collection of 4 or less lines in general position in
$\bbP^2$. From the definitions in section 1 it follows easily that,
$C$ is geometrically rigid. Moreover, it is equally clear that
collections of five or more lines in general position in $\bbP^2$ are
{\em not} geometrically rigid. Geometrically rigid divisors in
$\bbP^2$ (and hence their singular loci) are defined over $\Qbar$ (see
lemma~\ref{groth}):
\begin{thm}\label{uniq}
  Let $C$ be any divisor in $\bbP^2$ defined over $\bbC$ which is
  geometrically rigid. There is an automorphism $g$ of $\bbP^2$ so
  that $g(C)$ is defined over $\Qbar$.
\end{thm}
Now, if $C$ is a curve of degree 1 or 2 in $\bbP^2$, then $C$ is
geometrically rigid but a general curve of degree 3 or more is not.
In spite of this we will see that there are many geometrically rigid
divisors in $\bbP^2$. In fact (see the end of Section 3),
\begin{thm}\label{exist}
  Let $C$ be any divisor in $\bbP^2$ defined over $\Qbar$, and $T$ be
  a finite set of points in $\bbP^2$ defined over $\Qbar$. There is a
  geometrically rigid divisor $D$ in $\bbP^2$ so that $C\subset D$ and
  $T$ is contained in the singular locus of $D$.
\end{thm}
These results give a geometric characterisation of reduced algebraic
subschemes of $\bbP^2$ that are defined over $\Qbar$. As an easy
corollary we have a generalisation of Belyi's characterisation to the
case of surfaces.
\begin{cor}
  Let $S$ be a smooth projective surface, $C$ a divisor in $S$ and $T$
  a finite set of points in $S$.
  
  Assume that $S$, $C$ and $T$ are defined over $\Qbar$, then there is
  a geometrically rigid divisor $D$ in $\bbP^2$ and a finite morphism
  $f:S\to\bbP^2$ so that the image of $C$ and the branch locus of $f$
  are contained in $D$; moreover, the image of $T$ is contained in the
  singular locus of $D$.
  
  Conversely, suppose there is a tuple $(S,C,T)$ as above over $\bbC$
  and a finite morphism $f:S\to\bbP^2$ so that the image of $C$ and
  the branch locus of $f$ are sub-divisors of a geometrically rigid
  divisor $D$ and the image of $T$ is contained in the singular locus
  of $D$.  Then the tuple $(S,C,T)$ is isomorphic to (the base-change
  to $\bbC$ of) a tuple $(S_0,C_0,T_0)$ which is defined over $\Qbar$.
\end{cor}
It is reasonably clear that these results should be extendable {\em
mutatis mutandi} to higher dimension.
\begin{ack*}
  These results emerged during a seminar discussion with Gautham
  Dayal, Madhav Nori and G.~V.~Ravindra. I thank them for their
  valuable comments and criticisms. N.~Mohan Kumar made some valuable
  criticisms regarding section 2 which led me to look at the papers of
  Zariski more closely. N.~Fakruddin suggested lemma~\ref{nonempty}
  and the consequent simplification of lemma~\ref{groth}.
\end{ack*}

\section{Geometric Rigidity}
\begin{quote}
  Throughout the paper we work with schemes of finite type over a
  field of characteristic zero.
\end{quote}

Let $A$ be a smooth family of divisors in a smooth surface $S$; in
other words let $\calC\subset \calS=A\times S$ be a divisor with $A$
smooth. More generally, we can consider the case of non-constant
ambient spaces by assuming only that $\calS\to A$ is a smooth
projective morphism. We are interested in {\em topologically trivial}
families $p:(\calS,\calC)\to A$.  Over the field of complex numbers
this can be characterised by saying that any point $a\in A$ has an
analytic neighbourhood $U$ so that the pair $(U\times S,p^{-1}U)$ is
homeomorphic over $U$ to $U\times(S,p^{-1}(a))$. The geometric notion
of {\em equisingular} families results in topologically trivial
families.
\begin{rmk}
  The notion of equisingularity was first defined and studied by
  Zariski in a series of papers~\cite{equiI,equiII}. Theorem~7.4
  in~\cite{equiII} proves the equivalence of his definition with that
  studied here. Alternatively, one can directly prove Lemmas\ 
  \ref{nonempty},\ \ref{component} and\ \ref{groth} using his
  definitions. We require a specialised application of Zariski's
  results which we develop in this section.
\end{rmk}
A special case is that of a {\em
  family of divisors with normal crossings} which is characterised by
the following properties.
\begin{enumerate}
  \item The divisor $\calC$ is a divisor with normal crossings in
    $\calS$. 
  \item Each component of $\calC$ is smooth over $A$.
  \item The critical locus of $\calC\to A$ is \'etale over $A$.
\end{enumerate}
In particular, each component of the critical locus of $\calC\to A$
meets and is contained in exactly two components of $\calC$.

Now let $\calS_n\to\dots\to\calS_0=\calS$ be a sequence of blow-ups
with irreducible reduced centres $A_k\subset\calS_k$ such that $A_k\to
A$ is finite \'etale. Moreover, let $\calC_k$ denote the reduced union
of the total transform of $\calC$ in $\calS_k$ and the exceptional
locus of $\calS_k\to\calS$. We further assume that either,
\begin{enumerate}
\item $A_k$ is contained in the critical locus of $\calC_k\to A$ or,
\item $A_k$ is contained in $\calC_k$ but misses the critical locus of
  $\calC_k\to A$ entirely or,
\item $A_k$ lies in the complement of $\calC_k$ in $\calS_k$.
\end{enumerate}
While the latter two are irrelevant to the desingularisation it is
useful to allow these to simplify the proofs.  If $\calC_n$ is family
of divisors with normal crossings, then we call such a sequence of
blow-ups a {\em simultaneous desingularisation} of the family of
divisors $\calC\to S$. If such a sequence of blow-ups exists then we
say that the family is {\em simultaneously desingularisable} or {\em
  equisingular}. In order to understand how one arrives at this
definition we state
\begin{lemma}\label{nonempty}
  Fix a ground field $k$ of characteristic zero. Let $\calS\to A$ be a
  smooth family of projective surfaces of a reduced scheme $A$. Let
  $\calC\subset\calS$ be a reduced divisor. There is an open dense
  subset $U$ of $A$ over which $\calC$ is an equisingular family.
\end{lemma}
\begin{proof}
  We can replace $A$ by its smooth locus and further operate on each
  component individually; thus we can assume that $A$ is smooth and
  irreducible. Now, consider the reduced critical locus of $\calC\to
  A$. This is a closed subscheme $B$ of $\calC$ which is generically
  finite over $A$. Thus the locus where $B\to A$ is not \'etale is a
  proper closed subscheme of $A$. We can replace $A$ by the complement
  of this closed subscheme. Now we can take $A_1=B$ and perform a
  blow-up of $\calS$ along $A_1$ to obtain $\calS_1$. Since $A_1$ is
  \'etale over $A$ the resulting family $\calS_1\to A$ is smooth. Let
  $\calC_1$ denote the (reduced) union of the strict transform of
  $\calC$ in $\calS_1$ and the exceptional locus of the blow-up. We
  can now inductively construct the sequence $\calS_n$ as above. By
  the embedded desingularisation of curves in characteristic zero,
  there is an $n$ so that the generic fibre of $\calC_n\to A$ is a
  divisor with normal crossings; \ie each irreducible component (not
  geometrically irreducible component) of this generic fibre is smooth
  over the function field of $A$ and at most two of them meet at any
  singular point (which is closed over the function field of $A$) and
  this meeting is transversal. Now replace $A$ by the open subset
  where the critical locus of $\calC_n\to A$ is \'etale and each
  component of $\calC_n\to A$ is smooth. It follows that $\calC_n\to
  A$ is a family of divisors with normal crossings in $\calS_n\to A$.
\end{proof}
One point that is important from our perspective is the fact that $U$
is defined over $k$ since all schemes are of finite type over $k$. We
also note the following lemma.
\begin{lemma}\label{etale}
  Let $B_k$ be the image of the critical locus $B_n$ of $\calC_n\to A$
  in $\calS_k$ for each $k$. Then $B_n\to B_k$ and $B_k\to A$ are
  \'etale. Any component of $B_k$ that meets $A_k$ is actually $A_k$.
  Let $\calD_k$ be a union of components of $\calC_k$. If $\calD_k$
  and a component of $B_k$ meet then the latter is contained in the
  former. Finally, the critical locus of $\calD_k\to A$ is a union of
  components of $B_k$.
\end{lemma}
\begin{proof}
  We prove the statements by downward induction on $k$; we start at
  $k=n$ where this is true by the definition of a family of divisors
  with normal crossings. Now suppose that the result is proved for
  $B_{k+1}$ and for all divisors of the form $\calD_{k+1}$. Let $E_k$
  be the exceptional locus of $\calS_{k+1}\to\calS_k$. Then $E_k$ is
  contained in $\calC_{k+1}$ by the definition of $\calC_{k+1}$. The
  map $E_k\to A$ factors through $A_k\to A$ which is \'etale.
  
  Let $Y$ be the union of those connected components of $B_k$ which
  meet $A_k$; in particular, this includes those components which
  contain points where $B_{k+1}\to B_k$ is not an isomorphism. Let $X$
  be the inverse image of $Y$ in $B_{k+1}$; by the induction
  hypothesis $X\to A$ is \'etale. Moreover, each component of $X$
  meets $E_k$. By choosing $\calD_{k+1}=E_k$ we see that $X$ is
  contained in $E_k$ by the induction hypothesis. Thus, the morphism
  $X\to A$ is \'etale and factors through $A_k\to A$. It follows that
  $Y=A_k$.  Thus $B_k$ is the disjoint union of $A_k$ and components
  disjoint from $A_k$. The remaining components descend isomorphically
  from components of $B_{k+1}$ and $B_{k+1}\to A$ is \'etale by
  induction. Hence $B_k\to A$ is \'etale.
  
  Let $\calD_k$ be a union of irreducible components of $\calC_k$ and
  suppose that $\calD_k$ meets $A_k$. Let $\calD_{k+1}$ be its strict
  transform in $\calS_{k+1}$. Then $\calD_{k+1}$ must meet $E_k$; let
  $Z$ be any component of $\calD_{k+1}\cap E_k$. This is a divisor in
  $E_k$ which is contained in the critical locus of $\calD_{k+1}\cup
  E_k\to A$. By the induction hypothesis applied to $\calD_{k+1}\cup
  E_k$ we see that $Z$ is a component of $B_{k+1}$.  Hence, $Z\to A$
  is \'etale by induction, and the image of $Z$ is $A_k$ as above.
  Thus $\calD_k$ contains $A_k$.
  
  Finally, any critical point $p$ of $\calD_k\to A$ which is not the
  image of a critical point of $\calD_{k+1}\to A$, would have to lie
  in $A_k$. Either (a) there are two points $q$ and $q'$ that lie in
  $\calD_{k+1}\cap E_k$ over $p$, or (b) there is a point $q$ in
  $\calD_{k+1}\cap E_k$ where this interesection is not transvesal.
  In case (a) let $Z$ and $Z'$ be the components of $\calD_{k+1}\cap
  E_k$ that contain $q$ and $q'$ respectively ($Z=Z'$ is a
  possibility). Then $Z\to A_k$ and $Z'\to A_k$ are \'etale as
  explained above. In particular, $\calD_k\to A$ has critical points
  along $A_k$. In case (b), let $Z$ be the component of
  $\calD_{k+1}\cap E_k$ that contains $q$. The map $Z\to A_k$ is
  \'etale as above, hence $Z$ is smooth. Thus the intersection of
  $\calD_{k+1}$ and $E_k$ is non-transversal everywhere along $Z$.
  Thus, in this case $A_k$ is contained in the critical locus of
  $\calD_k\to A$ again. Any critical point of $\calD_k\to A$ is is
  thus either contained in $A_k$ which is contained in this critical
  locus or contained in the image of the critical locus of
  $\calD_{k+1}\to A$ which is a union of components of $B_k$. Since
  $A_k$ is contained in $B_k$ in both cases (a) and (b), it follows
  that the critical locus of $\calD_k\to A$ is a union of components
  of $B_k$.
\end{proof}
In particular, note that this means that $A_k$ is a connected
component of the critical locus of $\calC_k\to A$ if it meets this
locus; this strengthens the condition (1) in the definition above. The
fundamental lemma that we will use in our constructions is a corollary
of the above lemma.
\begin{lemma}\label{component}
  Let $(\calS,\calC)\to A$ be an equisingular family of divisors in a
  family of smooth projective surfaces over a smooth variety $A$. Let
  $\calD\subset\calC$ be a union of components of $\calC$, then
  $(\calS,\calD)\to A$ is an equisingular family of divisors.
\end{lemma}
\begin{proof}
  Let $\calS_n\to\dots\to\calS_0=\calS$ be a simultaneous
  desingularisation of $\calC$ as above. Let $\calD_k$ be the reduced
  total transform of $\calD$ in $\calS_k$. Since $\calD_n$ is a union
  of components of $\calC_n$, it too is a relative divisor with normal
  crossings over $A$. By above lemma we see that whenever $\calD_k\to
  A$ has a critical point on $A_k$, then $A_k$ is contained in this
  critical locus. Moreover, if $\calD_k$ meets $A_k$ then it contains
  it. Thus the given sequence of blow-ups is a simultaneous
  desingularisation of $\calD_k$.
\end{proof}

Let $C\subset S$ be a divisor. Let $G$ be an algebraic group of
automorphisms of $S$. Given a morphism $A\to G$, we can construct an
equisingular family containing $C$ as follows. Let $m:A\times S\to S$
denote the action of $A$ on $S$ and let $\calC=m^{-1}(C)$. More
generally, we say that a family $\calC\subset A\times S$ is {\em $G$
  iso-trivial}, if it is associated with a $G$-torsor on $S$. In other
words, each point $a\in A$ has an \'etale neighbourhood $B\to A$ so
that $\calC_B=\calC\times_A B$ is isomorphic over $B$ to
$m_B^{-1}(C_a)$ for some morphism $m_B:B\to G$. Any iso-trivial family
is clearly equisingular.

We now define $C$ to be a {\em geometrically rigid} divisor in $S$ if
this is the only way to construct equisingular deformations of $C$;
\ie for any equisingular family $\calC\subset A\times S$ parametrised
by a smooth connected variety $A$ so that $C$ is the fibre $p^{-1}(a)$
for some point $a$ in $A$, there is an algebraic group $G$ of
automorphisms of $S$ so that the family $\calC\to A$ is $G$ iso-trivial.

The following lemma follows easily from the construction of
universal deformations of divisors and the flattening stratification.
\begin{lemma}\label{groth}
  Let $S$ be smooth surface over an algebraically closed field $k$ and
  $C$ be a geometrically rigid divisor in $S$ defined over an
  algebraically closed extension $K$ of $k$. Then there is an
  automorphism $g$ of $S$ over $K$, so that $g(C)$ is the base change
  to $K$ of a curve $C_0$ in $S$ which is defined over $k$.
\end{lemma}
As a consequence, geometric rigidity is a sufficient criterion to
reduce the field of definition.
\begin{proof}
  Let $H$ be the Hilbert scheme of divisors in $S$ over $k$. Let $A$
  be the closure of the (non-closed) point of $H$ which corresponds to
  $C$. Then $A$ is a scheme of finite type over $k$ to which we can
  apply lemma~\ref{nonempty} above. Thus replacing $A$ an an open
  subscheme $U$ defined over $k$ we have an equisingular family
  $\calC\to A$ in $S\times A$ with generic fibre isomorphic to the
  given $C$.
  
  By the geometric rigidity of $S$ it follows that this family is
  isotrivial for some algebraic group $G$ of automorphisms of $S$.
  Thus there is a finite \'etale cover $A'\to A$ so that the family is
  group-theoretically trivial over $A'$. Since $k$ is algebraically
  closed there is a $k$-valued point of $A'$. The fibre of $\calC$ at
  this point is then a ``model'' of $(S,C)$ which is defined over $k$.
\end{proof}
In particular, we note that Theorem~\ref{uniq} follows.

\section{Constructions}

We now give inductive constructions of geometrically rigid
divisors to prove Theorem~\ref{exist}.
\begin{lemma}\label{addaline}
  Let $D$ be a geometrically rigid divisor in $\bbP^2$ and let $p$,
  $q$ be singular points of $D$. The divisor $D\cup\overline{pq}$ is
  geometrically rigid, where $\overline{pq}$ is the line joining the
  points $p$ and $q$.
\end{lemma}
\begin{proof}
  Let $\calC\to A$ be an equisingular deformation of
  $D\cup\overline{pq}$. We wish to construct a group-theoretic
  trivialisation of this deformation over a finite \'etale cover of
  $A$. 
  
  Let $A_1\to A$ (respectively $A_2\to A$) be a component of the
  critical locus of $\calC\to A$ which contains $p$ (respectively
  contains $q$). These are \'etale covers of $A$ by lemma~\ref{etale}.
  Let $B\to A$ be a connected \'etale cover of $A$ that dominates both
  covers; we have natural morphisms $P:B\to\bbP^2$ and $Q:B\to\bbP^2$
  passing through $p$ and $q$ respectively. Let $\calL\to B$ be the
  component of $\calC_B=\calC\times_A B$, that contains
  $\overline{pq}$. Then, the fibre of $L$ over $b\in B$ consists of
  the line joining $P(b)$ and $Q(b)$. Let $\calD_B$ be the union of
  the remaining components of $\calC_B$. By lemma~\ref{component}, the
  familiy $\calD_B\to B$ is an equisingular deformation of $D$.
  
  Now, by the geometric rigidity of $D$, we see that $\calD_B\to B$ is
  iso-trivial. In particular, we take a further \'etale cover (which
  we also denote by $B$ by abuse of notation) so that the family
  $\calD_B$ is group-theoretic. Now, $P(B)$ and $Q(B)$ continue to be
  part of the critical locus of $\calD_B\to B$, thus by the
  connected-ness of $B$ the trivialisation of the family must take
  them to $B\times \{p\}$ and $B\times \{q\}$ respectively. But then
  the same trivialisation also takes $L$ to the $B\times
  \overline{pq}$. Thus we have a group-theoretic trivialisation of
  $\calC_B$.
\end{proof}
Starting with the geometrically rigid divisor $Q$ of 4 lines in
general position on $\bbP^2$, we look at all the divisors obtained by
iterated application of the above lemma. The usual constructions of
projective geometry that give the field operations for points on a
line give the following result.
\begin{propn}\label{rational}
  Let $T$ be any finite set of points in $\bbP^2$ defined over
  $\bbQ$. There is a geometrically rigid divisor $D$ consisting of
  lines so that $T$ is contained in the singular locus of $D$.
\end{propn}
\begin{proof}
  Fixing the reference quadrilateral $Q$ consisting of four general
  lines in $\bbP^2$ also fixes a coordinate system so that the lines
  are given by $X=0$, $Y=0$, $Z=0$ and $X+Y+Z=0$. The singular points
  of the quarilateral are $(1:0:0)$, $(0:1:0)$, $(0:0:1)$, $(1:-1:0)$,
  $(1:0:-1)$ and $(0:1:-1)$.
  
  For any $t\in\bbP^2(\bbQ)$ and a geometrically rigid divisor $C_0$
  containing $Q$ we will construct a larger geometrically rigid
  divisor that contains $t$. We can then construct $D$ by starting
  with $Q$ and succesively adding each point of the finite set $T$.
  
  Thus we can assume that $T$ consists of just one point. Since at
  least one coordinate of $t$ is non-zero we can assume that it takes
  the form $(u:v:1)$ in these coordinates for some rational numbers
  $u$ and $v$.
  
  Now, suppose that we can add to $C_0$ and produce a geometrically
  rigid divisor $C$ so that the singular locus of $C$ contains
  $(u:0:1)$ and $(0:v:1)$. We can then add to $C$ the line $L$ joining
  $(u:0:1)$ and $(0:1:0)$, and the line $M$ joining $(0:v:1)$ and
  $(1:0:0)$, again producing a geometrically rigid divisor $C\cup
  L\cup M$ by the lemma~\ref{addaline}. Now the point $t$ is the
  intersection point of $L$ and $M$ so it is a singular point of this
  divisor as required.
  
  Similarly, if we can add to $C_0$ to produce a geometrically rigid
  divisor $C$ containing $(v:0:1)$ in its singular locus then the
  divisor $C\cup L$ is also geometrically rigid, where $L$ is the line
  joining $(v:0:1)$ and $(1:-1:0)$. The point $(0:v:1)$ whcih is the
  point of intersection of $L$ and the line $X=0$, is a singular point
  of this divisor. Thus to prove the result, it is enough to construct
  for each rational number $u$ a divisor $C_u$ containing $C_0$ so
  that the point $(u:0:1)$ is in the singular locus of $C_u$.
  
  We write $u=p/q$ where $q$ is a positive integer and $p$ is some
  integer. Suppose that we can construct a divisor $C$ containing $C_0$
  so that $(0:p:1)$ and $(0:-q:1)$ are singular points of $C$. Let $L$
  be the line joining $(1:0:-1)$ and $(0:-q:1)$; as before $C\cup L$
  is a geometrically rigid divisor. Moreover, $(1:-q:0)$ is a
  singular point of this divisor at it lies on $L$ and the line
  $Z=0$. Let $M$ be the line joining $(0:p:1)$ and $(1:-q:0)$; as
  before the divisor $C\cup L\cup M$ is geometrically rigid. The point
  $(p/q:0:1)$ is a singular point of this divisor as it lies on $M$
  and the line $Y=0$.
  
  Thus we have finally reduced to the problem of constructing for each
  integer $p$ a geometrically rigid divisor $C_p$ containing $C_0$ for
  which $(0:p:1)$ is a singular point. We will do this by induction on
  the absolute value of $p$. Let $L_1$ be the line joining $(0:1:0)$
  and $(-1:0:1)$, let $L_2$ be the line joining $(-1:1:0)$ and
  $(0:0:1)$. By the lemma~\ref{addaline} the divisor $Q\cup L_1\cup
  L_2$ is geometrically rigid. The point $(-1:1:1)$ is the
  intersection point of $L_1$ and $L_2$, hence it is a singular point
  of this divisor. Let $M$ be the line joining this point to
  $(1:0:0)$. Then $Q\cup L_1\cup L_2\cup M$ is geometrically
  rigid. The point $(0:1:1)$ is the intersection point of $M$ and the
  line $X=0$. Thus we have produced $C_p$ for every $p$ less than $1$
  in absolute value.

  Now, suppose that we have constructed a divisor $C$ containing
  $C_0$ which has $(0:p:1)$ and $(0:q:1)$ in its singular locus. Let
  $L_1$ be the line joining $(0:p:1)$ with $(1:0:0)$ and $L_2$ be the
  line joining $(-1:0:1)$ with $(0:1:0)$. The point $(-1:p:1)$ is a
  singular point of the geometrically rigid divisor $C\cup L_1\cup
  L_2$. Let $M_1$ be the line joining $(0:0:1)$ and $(-1:p:1)$r;
  the point $(-1:p:0)$ is a singular point of the geometrically rigid
  divisor $C\cup L_1\cup L_2\cup M_1$. Let $M_2$ be the line joining
  $(-1:p:0)$ and $(0:q:1)$; then $(-1:p+q:1)$ is a singular point of
  the geometrically rigid divisor $C\cup L_1\cup L_2\cup M_1\cup
  M_2$. Finally, we add the line $M_3$ joining $(1:0:0)$ and
  $(-1:p+q:1)$ to obtain a geometrically rigid divisor which has
  $(0:p+q:1)$ as a singular point.
  
  For every $p>1$ we apply the latter construction to the divisor
  $C_{p-1}\cup C_{1}$ each of which we have already constructed by the
  induction hypothesis and has singular points at $(0:p-1:1)$ and
  $(0:1:1)$. For $p<1$ we apply the latter construction to $C_{p+1}$
  which has $(0:p+1:1)$ and $(0:-1:1)$ as singular points. This
  provides the required construction and hence the result is proved.
\end{proof}

To construct points with coordinates in algebraic number fields we
need to have curves of degree greater than one in geometrically rigid
divisors.
\begin{lemma}\label{curve}
  Let $D$ be a geometrically rigid divisor on a rational surface $S$
  and let $T$ be a finite subset of the singular points of $D$. Let
  $L$ be a divisor class on $S$ so that the linear system $\abs{L-T}$
  has a unique element $E$. Then the divisor $D\cup E$ is
  geometrically rigid.
\end{lemma}
Actually, we only need the regularity of $S$ (\ie
$\cohom^1(S,\Oh_S)=0$) in the proof given below. Further
generalisations even for irregular surfaces are possible.
\begin{proof}
  Let $\calC\to A$ be an equisingular deformation of the divisor
  $C=D\cup E$. Let $B$ be the connected component of the critical
  locus of $\calC\to A$ that contains $T$. This is finite \'etale over
  $A$ by lemma~\ref{etale}. By base change we may assume that $B\to A$
  is an isomorphism. Thus, we can write $\calC=\calD\cup\calE$ where
  $\calD$ is the union of irreducible components of $\calC$ that meet
  $D$ and $\calE$ is the union of the irreducible components of
  $\calC$ that meet $E$. By the lemma~\ref{component}, $\calD\to A$ is
  an equisingular deformation of $D$. Thus by base change we have a
  group-theoretic trivialisation of $\calD$. Since $B$ is contained in
  the critical locus of $\calD\to A$, it is mapped into $T$ by the
  trivialisation. Thus, after applying this trivialisation, $\calE\to
  A$ becomes a family of divisors containing $T$.

  Now, the divisor class $L$ has no deformation since $S$ is
  rational. Thus, the divisor class of every fibre of $\calE\to A$ is
  in the class $L$. By assumption, $E$ is the unique such class
  containing $T$, thus $\calE\to A$ is the trivial family. Hence the
  trivialisation for $\calD\to A$ in fact gives a trivialisation of
  $\calE$ and $\calC$ as well.
\end{proof}

The above lemma allows us to apply the Lagrange interpolation formula to
prove the following proposition.
\begin{propn}
  Let $T$ be a finite set of algebraic points on $\bbP^2$, then there
  is a geometrically rigid divisor $D$ so that $T$ is contained in the
  singular locus of $D$.
\end{propn}
\begin{proof}
  As in the proof of proposition~\ref{rational}, given a geometrically
  rigid divisor $C_0$ which contains the reference quadrilateral $Q$
  and a point $t\in\bbP^2(\Qbar)$, we construct a larger divisor
  $C\contains C_0$ so that $t$ is in the singular locus of $t$. Since
  $T$ is a finite set we can inductively add all the points $t\in T$
  to obtain the required divisor $D$.  Thus we can assume that $T$
  consists of one point $t$.
  
  Again, as in the proof of proposition~\ref{rational} we can further
  reduce to the case where the point has the form $(u:0:1)$ where $u$
  is an algebraic number. Let $f(T)$ be a monic polynomial with
  rational coefficients for which $f(u)=0$; let $n$ be the degree of
  $f$. Let $F$ be the set of points $(k:f(k):1)$ for $k=0,\ldots,n^2$.
  The curve $E$ defined by $YZ^{n-1}=f(X/Z)Z^n$ passes through these
  $n^2+1$ points. Thus it is the unique curve of degree $n$ that
  passes through these points. Let $C$ be a divisor (containing the
  quadrilateral $Q$) constructed using proposition~\ref{rational}
  which contains $F$ in its singular locus.  The lemma~\ref{curve}
  then asserts that $D=C\cup E$ is geometrically rigid. The point
  $(u:0:1)$ is a point of intersection of $E$ and the line $Y=0$ which
  lies in $Q$; hence it is a singular point of $D$.
\end{proof}
Finally, any curve of degree $n$ defined over $\Qbar$ is uniquely
determined in its divisor class by $n^2+1$ distinct $\Qbar$-valued
points on it.
\begin{proof} (of the Theorem~\ref{exist}). Let $C$ be any curve of degree $n$
  in $\bbP^2$ which is defined over $\Qbar$. Let $T$ be a collection
  of $n^2+1$ distinct points on this curve over $\Qbar$. Let $D$ be a
  geometrically rigid divisor in $\bbP^2$ that contains $T$ in its
  singular locus. By lemma~\ref{curve} the divisor $D\cup C$ is
  geometrically rigid. Applying this argument to each component of a
  given divisor in $\bbP^2$ defined over $\Qbar$, we have the 
  result. 
\end{proof}

\section{Remarks and Open Problems}
A similar collection of arguments can be used to obtain geometrically
rigid configurations in $\bbP^n$ for $n\geq 3$. Projection arguments
can be used to define the notion of equisingular deformations of
in higher (co-)dimensions.  Arguments similar to the ones in the
previous section can then probably be used to show:
\begin{prob}
  For each $k$ between $0$ and $n-1$, let $T_k$ be a closed subscheme of
  $\bbP^n$ of pure dimension $k$ that is defined over $\Qbar$. Then
  there is a geometrically rigid divisor $S_{n-1}$ in $\bbP^n$ so that
  if $S_k$ is defined inductively as the singular locus of $S_{k+1}$,
  then $S_k$ has pure dimension $k$ and $T_k\subset S_k$.
\end{prob}
Another possible generalisation of Belyi's theorem is the following:
\begin{prob}
  If $C$ is a projective algebraic curve over a field of transcendence
  degree $r$ is there a morphism $f:C\to\bbP^1$ for which the branch
  locus has cardinality less than or equal to $3+r$.
\end{prob}
Belyi's original arguments can be used to show that the branch locus
can be assumed to be defined over the field of rational functions in
$r$ variables. However, there does not seem to be any obvious way to
reduce the number of points to $3+r$. The converse (that such a cover
is defined over a field of transcendence degree at most $r$) follows
from the the fact that $s$-tuples of points in $\bbP^1$ have a moduli
space of dimension $s-3$.

Finally, it is clear from the above construction that the complexity
of the configuration required to obtain rigidity is related to the
height of the defining equation of a curve. Can this relation be
explicitly used to define a notion of height?

\providecommand{\bysame}{\leavevmode\hbox to3em{\hrulefill}\thinspace}

\end{document}